\documentclass[10pt]{article}

\usepackage{graphics}
\usepackage{graphicx}

\usepackage[a4paper, left=35mm,right=35mm,top=34mm,bottom=34mm]{geometry}
\usepackage[utf8]{inputenc}
\usepackage[T1]{fontenc}
\usepackage[english]{babel}

\usepackage{enumerate}
\usepackage{graphicx}
\usepackage{hyperref}
\hypersetup{
    colorlinks=true,
    linkcolor=blue,
    filecolor=magenta,      
    urlcolor=cyan,
}

\usepackage{mathtools,amsthm,amssymb,amsfonts}
\usepackage{algorithm}
\usepackage{algorithmic}
\makeatother
\theoremstyle{plain}

\theoremstyle{definition}

\theoremstyle{remark}

\usepackage{caption} 
\usepackage[tight,footnotesize]{subfigure}

\usepackage{fancyhdr}
\usepackage{array}
\usepackage[usenames, dvipsnames]{xcolor}
\usepackage{tikz}
\usepackage{pgfplots}
\usepackage{ecp/ecpmashpc}
\usepackage{ecp/matrix_tikz/mrgcommon}
\usepackage{ecp/matrix_tikz/matrix_subdomain}

\author{
  {\normalsize Abal-Kassim Cheik Ahamed}\thanks{CUDA Research Center, Applied Mathematics and Systems Laboratory, \'Ecole Centrale Paris, France.}
	\and
  {\normalsize Fr\'ed\'eric Magoul\`es}\thanks{CUDA Research Center, Applied Mathematics and Systems Laboratory, \'Ecole Centrale Paris, France
    (correspondence, frederic.magoules@hotmail.com).}
		}	
\title{Accelerated solution of Helmholtz equation with Iterative Krylov Methods on GPU}
\date{}

\begin{document}
\maketitle
\thispagestyle{fancy}

\begin{abstract}
\noindent  This paper gives an analysis and an evaluation of linear algebra operations on Graphics Processing Unit (GPU) with complex number arithmetics with double precision. Knowing the performance of these operations, iterative Krylov methods are considered to solve the acoustic problem efficiently. Numerical experiments carried out on a set of acoustic matrices arising from the modelisation of acoustic phenomena within a cylinder and a car compartment are exposed, exhibiting the performance, robustness and efficiency of our algorithms, with a ratio up to 27x for dot product, 10x for sparse matrix-vector product and solvers in complex double precision arithmetics.
\end{abstract}


\section{Introduction}
\label{sec:introduction}

Amongst other industries, acoustic performance is a major concern in automotive companies. To tackle this, several models are used. In this study, we limit ourselves to acoustic models applying to closed cavities where the acoustic problem is independent from the surrounding structure. We are then able to assume that the pressure field does not have any interactions with the enclosed structure.

Over the past decade, Graphics Processing Units (GPUs) have undergone the same rapid evolution as Central Processing Units (CPUs) in the past. Its extent is such as computers with high Floating point Operations Per Second (Flops) are now widely spread.
If when they appeared, the usage of GPUs were focused on computations associated with graphic display, they encountered a considerable evolution ten years later when General Purpose GPUs (GPGPUs) started to use graphic hardware for numerical computations. With the computing power of GPUs drastically increasing, science and engineering started to use GPUs massively. This was also due to NVIDIA releasing Compute Unified Device Architecture (CUDA)~\cite{GPU:CUDA4.0:2011}, which offered a high level GPGPU-based programming language.

Before the migration of GPGPU in 2000, most numerical simulations were carried out on CPU clusters. However, the release of CUDA enabled scientists to take advantage of the high computational power of GPUs: the CPU freed itself from heavy calculations by distributing them to the GPU. The workload became shared out as follows: the CPU managed the sequential blocks of the code while the GPU dealt with the parallel ones. This can be explained by their respective architecture. The simplified architecture of a CPU processor is composed of several memories with multiple levels of cache memories associated, a basic unit of computation and a more complex control unit.

In this paper, we use \emph{Alinea}~\cite{cheikahamed:2012:inproceedings-1,cheikahamed:2012:inproceedings-2}, our own research group library, which implements many algorithms in C++. \emph{Alinea} is a scalable library used both for linear algebra operations and more advanced operations, such as iterative Krylov on both CPU and GPU clusters for real and complex number arithmetics in single and double precision. In this paper we are concerned only with complex number arithmetics.
We propose an evaluation of the solution of sparse and large size linear systems with matrices, described by complex number arithmetics, and arising from the finite~\cite{magoules:journal-auth:7} or infinite~\cite{magoules:journal-auth:26,magoules:journal-auth:19,magoules:journal-auth:15} element discretisation of acoustic problems, modeled with Helmholtz equation. Iterative Krylov are well suited for this kind of problems. Nevertheless, these methods require the computation of many linear algebra operations, in particular sparse matrix-vector product (SpMV), which requieres high computational time.
To accelerate the solution of these sparse linear systems, a computationnal effort is required to perform linear algebra operations efficiently.

The remainder of the paper is organised as follows. Section~\ref{sec:background_motivation} briefly gives the background and motivation of this work. Section~\ref{sec:case_study_automotive_acoustic_flow} decribes the industrial test cases used to evaluate and analyze our work. Section~\ref{sec:basic_linear_algebra_operations} presents how complex number arithmetics are designed. In Section~\ref{sec:advanced_linear_algebra_operations} we give the evaluation of basic (addition of vectors, scale of vectors, etc.) and advanced (sparse matrix-vector multiplication, etc.) linear algebra operations required to perform iterative Krylov methods. Section~\ref{sec:matrix_vector_product} presents numerical experiments on iterative Krylov methods which clearly exhibit the robustness and effectiveness of our implementation on GPU for solving acoustic problem. Concluding remarks are given in Section~\ref{sec:concluding_remarks}

\section{Background and Motivation}
\label{sec:background_motivation}

If Graphics Processing Units (GPUs) were first used exclusively for graphics applications like Graphical User Interface, this has evolved in the last few decades and due to their heavy computational power, GPUs are now widely used to accelerate scientific computation by graphics card hardware. GPUs have become a major tool in scientific computation, thanks both to the rapid pace at which they have improved, and the flexibility offered by languages to program on them such as CUDA, an extension from C/C++. This is especially true when the time consumption of a numerical simulation is taken into account, and GPUs have therefore provided a huge boost to a multitude of applications of science~\cite{GPU:SC:2009} and engineering~\cite{creel_high_2012}.

Unlike Central Processing Units (CPUs), GPUs have an architecture with multiple arithmetic compute units, as well as several levels of memory. This enables them to compute simultaneously the same operation over a million times. For example, the last graphics card of the Kepler family, the K40, has a 4.29 teraflops single-precision and 1.43 teraflops double-precision peak floating point performance. The K40 also has a 12 GB storage memory, and more generally the storage memory of GPUs have also greatly increased with time. The execution time of a GPU algorithm is strongly impacted by both the configuration of the distribution of the threads on the grid~\cite{cheikahamed:2013:inproceedings-3} and the memory managment~\cite{ref:cheikahamed:2012:inproceedings-1:DARG:2013}.  For the most time consuming linear algebra operation, the sparse matrix vector product, the structure of the matrix and the type of the matrix storage format are both of most importance, as shown in~\cite{GPU:BG:2009,cheikahamed:2012:inproceedings-2}. This, along with the importance of the distribution of threads, is also proved in~\cite{cheikahamed:2012:inproceedings-1} on iterative Krylov methods on GPU used to solve linear systems with real number arithmetics.

The acoustic problem, in its simplest linear form, is governed in the frequency domain by the Helmholtz equation with suitable boundary conditions.
When high frequency regime is considered, the matrix of the linear system becomes very large.
The problem we are aimed to solve arises from the discretization of the Helmholtz equation in a bounded domain $\Omega$, with a boundary condition considered on the ouside boundary $\Gamma = \partial \Omega$. The Helmholtz equation is expressed as follows $-\nabla^2 u - k^2 u= g$, where $k = \frac{2\pi F}{c}$ is the wavenumber associated with the frequency $F\in\mathbb{R}$ and $c\in\mathbb{R}$ denotes the velocity of the medium that is different in space. In this paper, Dirichlet boundary conditions are considered along a part of $\Gamma$.
The frequency domain in which the solution is sought is usually limited, so as to analyze the acoustic response at specific places of the cavity (for instance around the driver's ears).
To carry this out, a suitable numerical model has to be used.
For complex geometries two models can be chosen, depending on the boundary conditions.
If there are conditions on all boundaries of the domain, then boundary element (BE) methods can be used. Else, finite element (FE) methods, which are methods based on the domain, are used to solve a weak formulation of the problem. When using FE methods, mesh requirements (around 10 nodes per wavelength are necessary) make the mesh sizes gigantic when dealing with high frequencies. This paper focuses on effectively handling large size acoustic problems using FE method.

\section{Case Study: automotive and cylinder acoustic}
\label{sec:case_study_automotive_acoustic_flow}

In this part, we present the finite element meshes used for solving the acoustic problems arising from the automotive industry~\cite{magoules:journal-auth:4}.

\begin{figure}[!ht]
\centering
\includegraphics[scale=0.14]{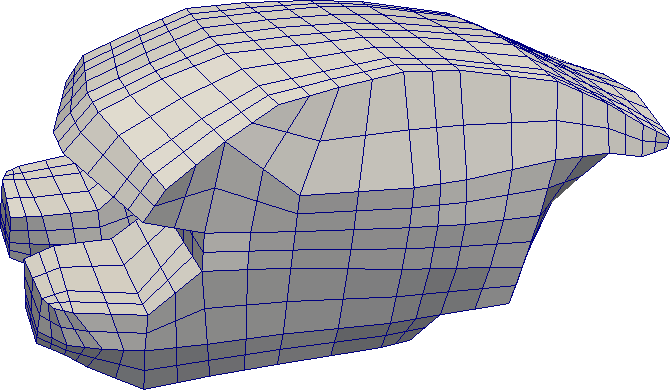}
\includegraphics[scale=0.14]{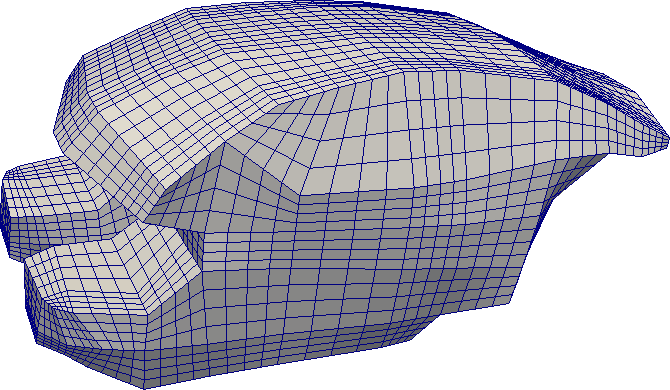}
\includegraphics[scale=0.14]{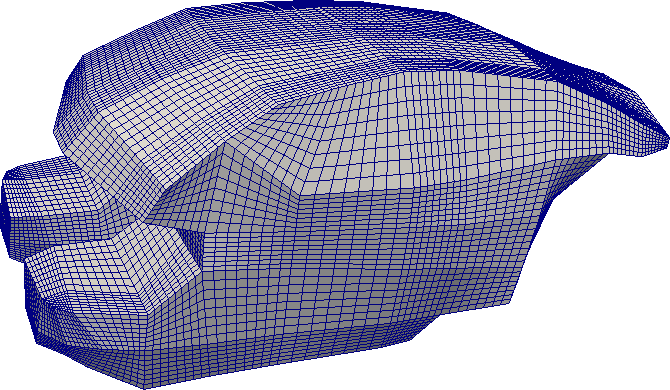}
\includegraphics[scale=0.14]{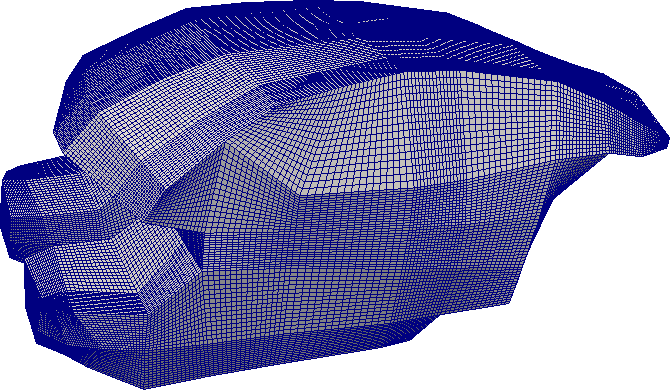}
\caption{Audi 3D, {\small h = (0.133425, 0.066604, 0.033289, 0.016643)}}
\label{fig:img:mesh_Audi3D}
\end{figure}

\begin{figure}[!ht]
\centering
\includegraphics[scale=0.17]{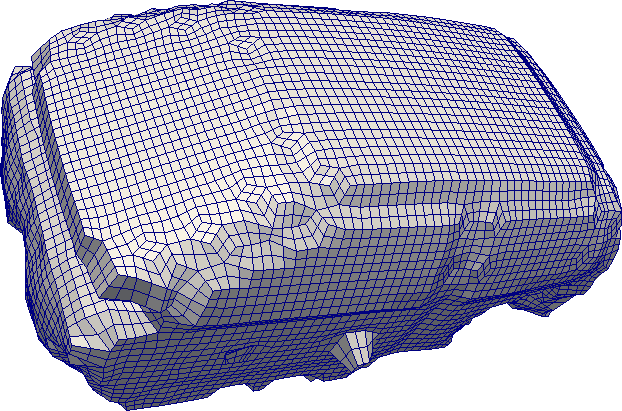}
\includegraphics[scale=0.17]{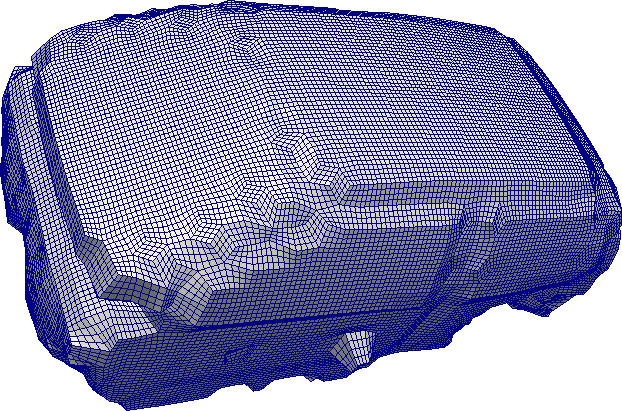}
\includegraphics[scale=0.17]{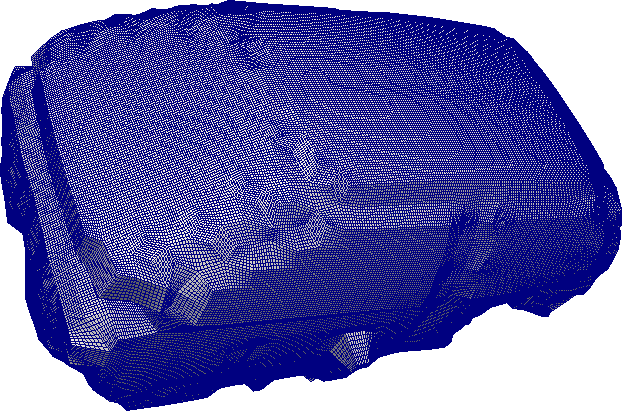}
\caption{Twingo 3D, {\small h = (0.077866, 0.038791, 0.019379)}}
\label{fig:img:mesh_Twingo3D}
\end{figure}

\begin{figure}[!ht]
\centering
\includegraphics[scale=0.2]{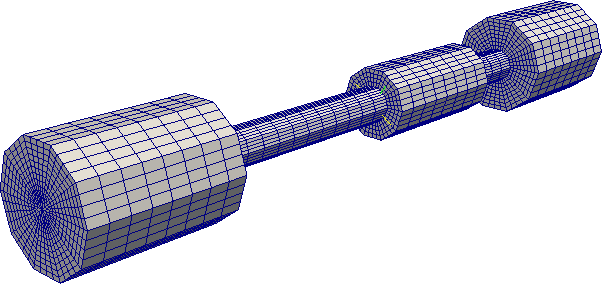}
\includegraphics[scale=0.2]{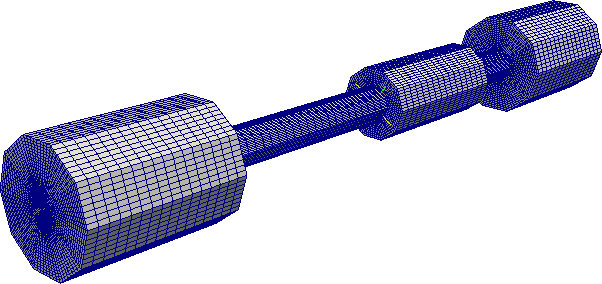}
\includegraphics[scale=0.2]{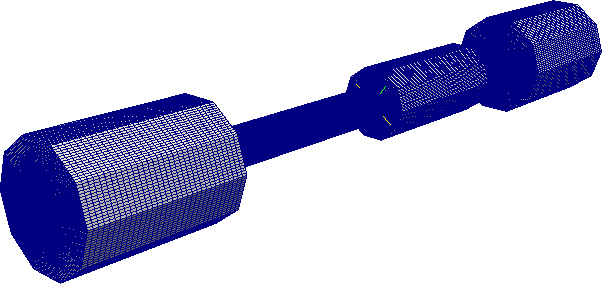}
\caption{Cylinder 3D, {\small h = (0.033973, 0.016949, 0.008342)}}
\label{fig:img:mesh_Cylinder3D}
\end{figure}

We now focus on two numerical examples that enable us to evaluate the performances of our procedures. We consider the study in a car compartment with Audi (Audi3D) and Twingo (Twingo3D), and in a cylinder (3D cylinder).
Let's look at the example of a \emph{car compartment}. The goal is to construct the frequency response function at driver's ear from the velocity boundary conditions along the firewall. Understanding this problem can help solve similar problems where the evaluation of the acoustic response to vibrating panels inside a cavity is at stake.
Several sources can explain such mechanical vibrations. The vibrations can indeed be air-borne or structural-borne. And the prediction of these vibrations can be a difficult task.
In the case of automotive applications, the higher the frequency is, the worst the quality of numerical predictions for mechanical vibrations is. Acoustic predictions depending on the treatment of these mechanical vibrations, precise acoustic predictions are possible only if correct vibration profiles along the car body are provided.
According to advanced FE methodologies used on car bodies, computing accurate results becomes difficult when the frequencies are higher than 2500 Hz.
Such difficulty to produce correct results at high frequencies can be explained by the complex mechanical structure of a car body. Usual models do not consider parameters which are essential to understand the behavior of a car body at high frequencies, such as the characteristics of the connections, the damping properties, etc. We can note that modifying the models to make them take such parameters into account is not easy.
At these high frequencies, some variability effects become important and complicate the predictions. For instance two car bodies produced in the same way, and that could be considered identical, may present drastically different vibro-acoustic behaviors at such frequencies.
The meshes of Audi 3D with different refinements (size $h$) are presented in Figure~\ref{fig:img:mesh_Audi3D}.
Figure~\ref{fig:img:mesh_Twingo3D} describes the meshes of Twingo 3D with different refinements 
(size $h$) .
In Figure~\ref{fig:img:mesh_Cylinder3D} are illustrated the Cylinder 3D meshes for three refinements (size $h$).

\subsection{Matrices tested}
\label{subsec:matrices_tested}

For the analysis proposed in this paper, we use matrices coming from the FE discretization of the Helmholtz equation for the study of car compartment and cylinder acoustic problems. These matrices are in stored Compressed Sparse Row (CSR).
Table~\ref{tab:sketches_matrices_audi_Twingo_Cynlinder_3D} collect a set of matrices respectively arising from Audi3D, Twingo3D and Cylinder3D. These tables give the properties of each considered matrix where $h$ is the size of the square matrix, $nz$ the number of nonzero elements, $density$ the density, i.e., the number of nonzero elements divided by the total number of matrix elements, $bandwidth$ the upper bandwidth which is equal to lower bandwidth for symmetric matrix, $max row$ the maximum row density, $nz/h$ the mean row density and $nz/h stddev$ the standard deviation of $nz/h$. The sparse matrix pattern and the coefficients distribution are respectively reported in the first column and second column.
\begin{table}
\centering
\scalebox{0.7}{
\begin{tabular}{m{6.4cm}r}
\hlinewd{1.0pt}
\infosparse{img/png/sparse/carHexOnlyFixed-GALERKIN-mat-1}{Audi3D-1}{1727}{16393}{0.550}{1436}{27}{9.492}{10.205}\\
\multicolumn{2}{c}{\emph{3D acoustic FE matrix. Audi car (mesh size = 0.133425, length wave = 3.5).}}\\
\hlinewd{1.0pt}
\infosparse{img/png/sparse/carHexOnlyFixed-GALERKIN-mat-2}{Audi3D-2}{11637}{188455}{0.139}{11237}{27}{16.194}{11.223}\\
\multicolumn{2}{c}{\emph{3D acoustic FE matrix. Audi car (mesh size = 0.066604, length wave = 3.5).}}\\
\hlinewd{1.0pt}
\infosparse{img/png/sparse/carHexOnlyFixed-GALERKIN-mat-3}{Audi3D-3}{85001}{1781707}{0.025}{84474}{27}{20.961}{9.832}\\
\multicolumn{2}{c}{\emph{3D acoustic FE matrix. Audi car (mesh size = 0.033289, length wave = 3.5).}}\\
\hlinewd{1.0pt}
\infosparse{img/png/sparse/carHexOnlyFixed-GALERKIN-mat-4}{Audi3D-4}{648849}{15444211}{0.004}{520461}{27}{23.802}{7.720}\\
\multicolumn{2}{c}{\emph{3D acoustic FE matrix. Audi car (mesh size = 0.016643, length wave = 3.5).}}\\
\hlinewd{1.0pt}
\infosparse{img/png/sparse/legacyCavity-GALERKIN-mat-0}{Twingo3D-0}{8439}{143889}{0.202}{6268}{27}{17.050}{11.047}\\
\multicolumn{2}{c}{\emph{3D acoustic FE matrix. Twingo car (mesh size = 0.077866, length wave = 9.5).}}\\
\hlinewd{1.0pt}
\infosparse{img/png/sparse/legacyCavity-GALERKIN-mat-1}{Twingo3D-1}{62357}{1351521}{0.035}{53935}{33}{21.674}{9.364}\\
\multicolumn{2}{c}{\emph{3D acoustic FE matrix. Twingo car (mesh size = 0.038791, length wave = 9.5).}}\\
\hlinewd{1.0pt}
\infosparse{img/png/sparse/legacyCavity-GALERKIN-mat-2}{Twingo3D-2}{479169}{11616477}{0.005}{470625}{39}{24.243}{7.233}\\
\multicolumn{2}{c}{\emph{3D acoustic FE matrix. Twingo car (mesh size = 0.019379, length wave = 9.5).}}\\
\hlinewd{1.0pt}
\infosparse{img/png/sparse/meshFig3-GALERKIN-mat-0}{Cylinder3D-0}{2717}{30969}{0.420}{2361}{75}{11.398}{11.453}\\
\multicolumn{2}{c}{\emph{3D acoustic FE matrix. Cylinder (mesh size = 0.033973, length wave = 9.5).}}\\
\hlinewd{1.0pt}
\infosparse{img/png/sparse/meshFig3-GALERKIN-mat-1}{Cylinder3D-1}{19041}{343677}{0.095}{18629}{75}{18.049}{11.051}\\
\multicolumn{2}{c}{\emph{3D acoustic FE matrix. Cylinder (mesh size = 0.016949, length wave = 9.5).}}\\
\hlinewd{1.0pt}
\infosparse{img/png/sparse/meshFig3-GALERKIN-mat-2}{Cylinder3D-2}{142049}{3151773}{0.016}{141289}{75}{22.188}{9.125}\\
\multicolumn{2}{c}{\emph{3D acoustic FE matrix. Cylinder (mesh size = 0.008342, length wave = 9.5).}}\\
\hlinewd{1.0pt}
\hlinewd{1.0pt}
\end{tabular}
}
\caption{Sketches of Audi, Twingo, Cylinder FE matrices}
\label{tab:sketches_matrices_audi_Twingo_Cynlinder_3D}
\end{table}

\subsection{Sparse matrix formats}
\label{subsec:sparse-matrix_formats}

\begin{figure}[!ht]
\centering
\includegraphics[scale=0.23]{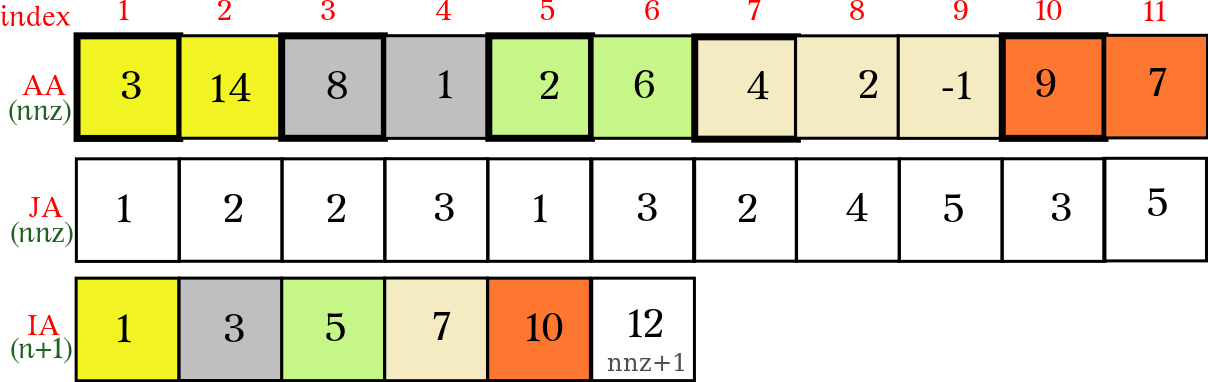}
\caption{Compressed sparse row storage format (CSR) of matrix, TABLE~\ref{tab:matrix_and_pattern}}
\label{fig:img:csr_format}
\end{figure}

The matrices of the acoustic problem are of large size and are sparse. Most of the elements are zero. 
The distribution of non-zero values depends on the properties of the initial problem.
To take the best advantage of the memory storage, sparse matrices are stored in compressed formats, i.e., only non-zero elements are allocated.
It exists divers storage formats~\cite{saad_iterative_2003} such as Compressed-Sparse Row (CSR)~\cite{GPU:BG:2009}, Coordinate (COO), ELLPACK (ELL), Hybrid (HYB), etc. 
In this paper, we have chosen CSR format to store matrices. The matrix is stored in three one-dimensional arrays, as drawn in Figure~\ref{fig:img:csr_format} and corresponds o the matrix shown Table~\ref{tab:matrix_and_pattern}. The first two arrays of size $nz$, $AA$ and $JA$ contain respectively the non-zero coefficients of the matrix in row major order and the column indices, hence $JA(j)$ corresponds to the column index in dense matrix $A$ of the coefficient $AA(j)$. The last array, $IA$, of size $n+1$, contains pointers to the start of each row. $IA(i)$ and $IA(i - 1) - 1$ match to the start and the end of the $i-th$ row in arrays $AA$ and $JA$, i.e., $IA(n + 1) = nz + 1$.
\begin{table}
\centering
\begin{tabular}{m{5cm}m{5cm}}
\scalebox{1.0}{
$A=\begin{pmatrix}
\colorbox{gray}{3} & \colorbox{gray}{14} & 0 & 0 & 0 \\
 0 & \colorbox{gray}{8} & \colorbox{gray}{1} & 0 & 0 \\
 \colorbox{gray}{2} & 0 & \colorbox{gray}{6} & 0 & 0 \\
 0 & \colorbox{gray}{4} & 0 & \colorbox{gray}{2} & \colorbox{gray}{-1} \\
 0 & 0 & \colorbox{gray}{9} & 0 & \colorbox{gray}{7}
\end{pmatrix}$
}
&
\scalebox{0.6}{
\SetDisplayColumnIndexD{1}{0}
\SetDisplayRowIndexD{1}{0}
\def\nnzcoef{1/1,1/2,2/2,2/3,3/1,3/3,4/2,4/4,4/5,5/3,5/5}
\DrawGivenMatrix{(0,0)}{5}{5}{\nnzcoef}
}\\
\end{tabular}
\caption{Left (matrix), Right (matrix pattern)}
\label{tab:matrix_and_pattern}
\end{table}

\subsection{Hardware configuration}
\label{subsec:experiments_configuration}

The experiments have been performed on workstation based on an Intel Core i7 920 2.67Ghz, which has 4 physical cores and 4 logical cores, 12GB RAM, and two nVidia graphics card: a Tesla K20c (device \#0) with 4799GB memory and GeForce GTX 570 with 1279MB memory (device \#1). The cards are double precision compatible. In this paper the devices Tesla K20c and GTX570 are respectively denoted \emph{gpu\#0} and \emph{gpu\#1}.

The fact that the accuracy of the native clocks of the GPU and the host are respectively a few nanoseconds and a few milliseconds can generate parasites on the measuring of execution time of our programs. To overcome it, we execute the same operations of our benchmark at least 10 times and at least enough times for the total measured time to exceed by at least 100 times the clock accuracy. For the sake of accuracy, we in fact run every operation 100 times, and the reported times correspond to the average time.

\section{Handling Complex Numbers on CUDA}
\label{sec:complex_number_on_CUDA}

The finite element discretization of the Helmholtz equations for acoustic problems leads to complex number arithmetics matrices.
GPU have originally been proposed for integer arithmetics. Most numerical simulations need floating point number artihmetics, which decrease the performance of computations.  
It is even worse in terms of floating point operations when we use numbers with double precision.
Since acoustic problems involve complex number arithmetics with double precision floating point operations, the awaited efficiencies are exceedingly low.
In order to get the best advantage of GPU architecture, we aim to best optimize the usage of complex number on CUDA. 
From the viewpoint of numerical calculation, a complex number is defined as a set of two real numbers, which corresponding to real and imaginary part. 
A natural way to represent a complex number in memory through the following structure
\begin{lstlisting}[language=c,basicstyle=\small,xleftmargin=25pt,caption={Complex number in memory},captionpos=b, label=lst:complex_number]
struct complex {
  double x; // real part
  double y; // imaginary part
};
\end{lstlisting}
Considering that the previous structure has not padding, both real numbers $x$ and $y$ are contiguous on memory, i.e., the offset of both is null. 
In our CUDA implementation, we keep the same design as previously. The native library of CUDA, \lstinline+libcudart.so+, offers \lstinline+cuDoubleComplex+ a \lstinline+double2+ structure (\lstinline+typedef double2 cuDoubleComplex+), which consists of the same procedure given in Listing~\ref{lst:complex_number}. 
For performance outcomes, we prefer to design a complex template class \lstinline+complex<T>+, which redefines all the operations available in standard \lstinline+std::complex+. The functions of the standard complex are called host, i.e., only executable by the CPU. To be executable on GPU, the functions are redefined with \lstinline+__host__ __device__+, which consists of both CPU and GPU code. We place our complex class template structure into a namespace \lstinline+stdmrg+ in order to avoid confusing it with the standard. This class can be used in both host and device platform. As in all CUDA code, in our library, there exists a copy function that allows to transfer data from host to device, and from device to host. This function offers the possibilities to copy from CPU with \lstinline+complex<T>+ or \lstinline+stdmrg::complex<T>+ to GPU with \lstinline+stdmrg::complex<T>+. Similarly, from GPU \lstinline+stdmrg::complex<T>+ to \lstinline+complex<T>+ or \lstinline+stdmrg::complex<T>+.

Distribution of threads is not an automated process. References~\cite{GPU:DZO:2011}~\cite{cheikahamed:2012:inproceedings-2} proved that the threading organization strongly impacts the performance of the numerical algorithm. GPU implementations using advanced gridification auto-tuning techniques, as developped in~\cite{cheikahamed:2012:inproceedings-1,cheikahamed:2012:inproceedings-2} are thus used in the following.

\section{Basic Linear Algebra Operations}
\label{sec:basic_linear_algebra_operations}

This section presents basic linear algebra algorithms, including assign of a vector, scale of vectors, element wise product, addition of vectors and dot product. We also collect the numerical experiments results we have performed to analyze the speed-up of the GPU code compared to the CPU code for complex number arithmetics with double precision. 

\begin{figure}[!ht]
\centering
  \begin{tikzpicture}[scale=0.75]
    \begin{axis}[
      height=8cm,
      width=10cm,
      xlabel=$size$,
      legend pos=north west,
      enlargelimits,
    ]
    \addplot[line width=2pt,dashed,black,mark=*] table[x index=0,y index=1,col sep=comma] {data/tikz/z_assign.txt};
    \addlegendentry{cpu time (ms)}
    \addplot[line width=2pt,red,mark=triangle] table[x index=0,y index=3,col sep=comma] {data/tikz/z_assign.txt};
    \addlegendentry{gpu\#0 time (ms)}
    \addplot[line width=2pt,blue,mark=*] table[x index=0,y index=5,col sep=comma] {data/tikz/z_assign.txt};
    \addlegendentry{gpu\#1 time (ms)}
    \end{axis}
  \end{tikzpicture}
  \begin{tikzpicture}[scale=0.75]
    \begin{axis}[
      height=8cm,
      width=10cm,
      xlabel=$size$,
      legend pos=south east,
      enlargelimits,
    ]
    \addplot[line width=2pt,dashed,black,mark=*] table[x index=0,y index=2,col sep=comma] {data/tikz/z_assign.txt};
    \addlegendentry{cpu GFlop/s}
    \addplot[line width=2pt,red,mark=triangle] table[x index=0,y index=4,col sep=comma] {data/tikz/z_assign.txt};
    \addlegendentry{gpu\#0 GFlop/s}
    \addplot[line width=2pt,blue,mark=*] table[x index=0,y index=6,col sep=comma] {data/tikz/z_assign.txt};
    \addlegendentry{gpu\#1 GFlop/s}
    \end{axis}
  \end{tikzpicture}

\caption{ZASSIGN [left: time in ms, right: GFlops]}
\label{fig:tikz:z_assign}
\end{figure}
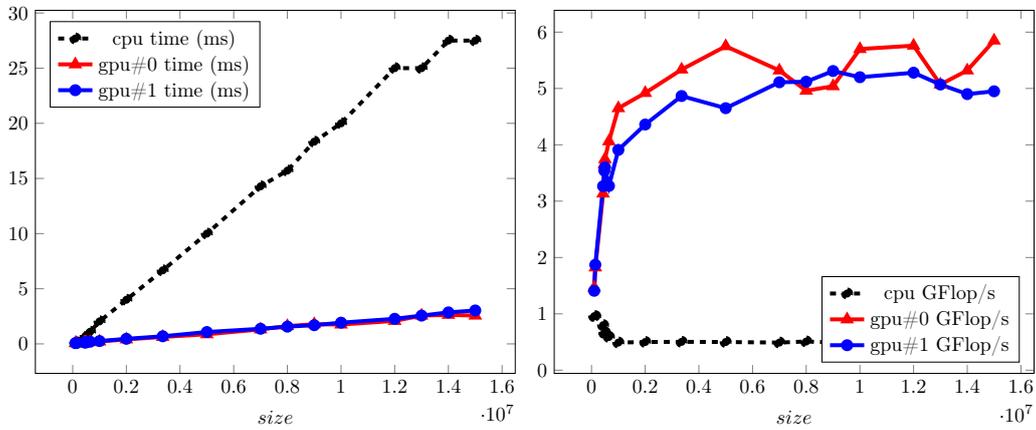

Table~\ref{tab:assign_tab} and Figure~\ref{fig:tikz:z_assign} give the complex double precision execution time in milliseconds (ms) of the \emph{assign operation}, where $h$ represents the size of the vector.
\begin{table}[!ht]  
\centering    
\renewcommand{\arraystretch}{1.3}
\renewcommand{\tabcolsep}{0.06cm}
\begin{tabular}{ccccccccc}
\hlinewd{1.0pt}
{\bf h} & {\bf cpu} & {\bf cpu} & {\bf gpu\#0} & {\bf gpu\#0} & {\bf gpu\#1} & {\bf gpu\#1} & {\bf ratio\#0} & {\bf ratio\#1} \\
{} & {\it time (ms)} & {\it Gflops} & {\it time (ms)} & {\it Gflops} & {\it time (ms)} & {\it Gflops} & {\it cpu/\#0} & {\it cpu/\#1} \\
\hlinewd{1.0pt}
100,000 & 0.10 & 0.96 & 0.07 & 1.47 & 0.07 & 1.41 & {\bf 1.53} & {\bf 1.47} \\
500,000 & 0.75 & 0.67 & 0.13 & 3.75 & 0.14 & 3.60 & {\bf 5.59} & {\bf 5.37} \\
1,000,000 & 2.04 & 0.49 & 0.22 & 4.65 & 0.26 & 3.91 & {\bf 9.49} & {\bf 7.98} \\
8,000,000 & 15.71 & 0.51 & 1.61 & 4.96 & 1.56 & 5.12 & {\bf 9.74} & {\bf 10.06} \\
10,000,000 & 20.00 & 0.50 & 1.75 & 5.70 & 1.92 & 5.20 & {\bf 11.40} & {\bf 10.40} \\
15,000,000 & 27.50 & 0.55 & 2.56 & 5.85 & 3.03 & 4.95 & {\bf 10.73} & {\bf 9.08} \\
\hlinewd{1.0pt}
\end{tabular}
\caption{Complex double precision Assign of vector (ZASSIGN)}
\label{tab:assign_tab}
\end{table}

The scale {\em scale operation} kernel is described as follows
\begin{lstlisting}[language=cuda,label=lst:xmy_2d]
__global__ void Scal(
        stdmrg::complex<double> alpha,
        const stdmrg::complex<double>* d_x, int size) {
  unsigned int x = blockIdx.x * blockDim.x + threadIdx.x;
  unsigned int y = threadIdx.y + blockIdx.y * blockDim.y;
  int pitch = blockDim.x * gridDim.x;
  int idx = x + y * pitch;
  if ( idx < size ) {
    d_x[idx] = alpha * d_x[idx];
  }
}
\end{lstlisting}
The running times of the scale {\em scale operation} for different size of vectors are given in Table~\ref{tab:scal_vectors_tab} and drawn in Figure~\ref{fig:tikz:z_scal}, where $h$ represents the size of the vectors.
\begin{table}[!ht]  
\centering    
\renewcommand{\arraystretch}{1.3}
\renewcommand{\tabcolsep}{0.06cm}
\begin{tabular}{ccccccccc}
\hlinewd{1.0pt}
{\bf h} & {\bf cpu} & {\bf cpu} & {\bf gpu\#0} & {\bf gpu\#0} & {\bf gpu\#1} & {\bf gpu\#1} & {\bf ratio\#0} & {\bf ratio\#1} \\
{} & {\it time (ms)} & {\it Gflops} & {\it time (ms)} & {\it Gflops} & {\it time (ms)} & {\it Gflops} & {\it cpu/\#0} & {\it cpu/\#1} \\
\hlinewd{1.0pt}
100,000 & 0.81 & 0.74 & 0.07 & 8.60 & 0.07 & 8.26 & {\bf 11.56} & {\bf 11.10} \\
500,000 & 4.17 & 0.72 & 0.15 & 19.50 & 0.17 & 17.16 & {\bf 27.08} & {\bf 23.83} \\
1,000,000 & 8.33 & 0.72 & 0.29 & 20.70 & 0.30 & 20.04 & {\bf 28.75} & {\bf 27.83} \\
8,000,000 & 65.00 & 0.74 & 1.75 & 27.36 & 2.22 & 21.60 & {\bf 37.05} & {\bf 29.25} \\
10,000,000 & 85.00 & 0.71 & 2.38 & 25.20 & 2.56 & 23.40 & {\bf 35.70} & {\bf 33.15} \\
15,000,000 & 120.00 & 0.75 & 3.13 & 28.80 & 3.70 & 24.30 & {\bf 38.40} & {\bf 32.40} \\
\hlinewd{1.0pt}
\end{tabular}
\caption{Complex double precision Scale of vectors (ZSCAL)}
\label{tab:scal_vectors_tab}
\end{table}
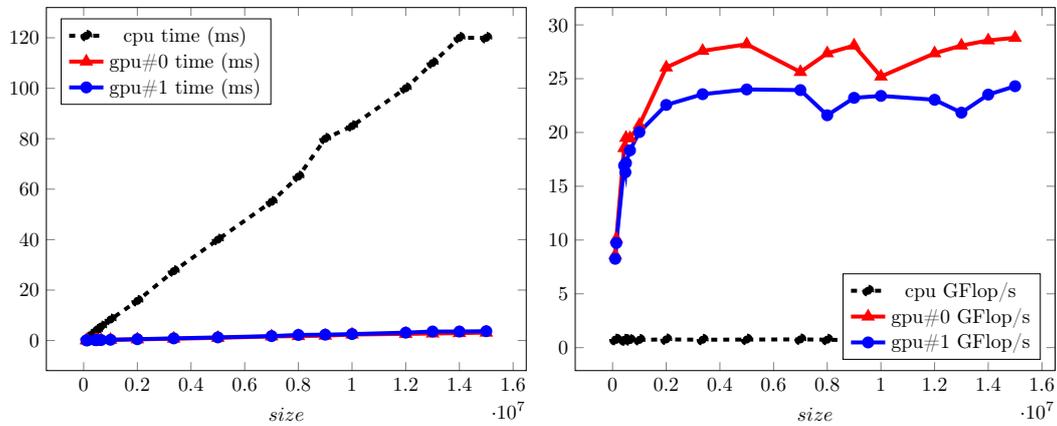
\begin{figure}[!ht]
\centering
  \begin{tikzpicture}[scale=0.75]
    \begin{axis}[
      height=8cm,
      width=10cm,
      xlabel=$size$,
      legend pos=north west,
      enlargelimits,
    ]
    \addplot[line width=2pt,dashed,black,mark=*] table[x index=0,y index=1,col sep=comma] {data/tikz/z_scal.txt};
    \addlegendentry{cpu time (ms)}
    \addplot[line width=2pt,red,mark=triangle] table[x index=0,y index=3,col sep=comma] {data/tikz/z_scal.txt};
    \addlegendentry{gpu\#0 time (ms)}
    \addplot[line width=2pt,blue,mark=*] table[x index=0,y index=5,col sep=comma] {data/tikz/z_scal.txt};
    \addlegendentry{gpu\#1 time (ms)}
    \end{axis}
  \end{tikzpicture}
  \begin{tikzpicture}[scale=0.75]
    \begin{axis}[
      height=8cm,
      width=10cm,
      xlabel=$size$,
      legend pos=south east,
      enlargelimits,
    ]
    \addplot[line width=2pt,dashed,black,mark=*] table[x index=0,y index=2,col sep=comma] {data/tikz/z_scal.txt};
    \addlegendentry{cpu GFlop/s}
    \addplot[line width=2pt,red,mark=triangle] table[x index=0,y index=4,col sep=comma] {data/tikz/z_scal.txt};
    \addlegendentry{gpu\#0 GFlop/s}
    \addplot[line width=2pt,blue,mark=*] table[x index=0,y index=6,col sep=comma] {data/tikz/z_scal.txt};
    \addlegendentry{gpu\#1 GFlop/s}
    \end{axis}
  \end{tikzpicture}

\caption{ZSCAL [left: time in ms, right: GFlops]}
\label{fig:tikz:z_scal}
\end{figure}

{\em Double-precision complex Alpha X Plus Y (Zaxpy)}, i.e., $y[i] = \alpha \times x[i] + y[i]$, is a level one (vector) operation between two complex number arithmetics vectors in the Basic Linear Algebra Subprograms (BLAS) package.
The addition of vectors is performed by a simple GPU kernel of the form:
\begin{lstlisting}[language=cuda,label=lst:zaxpy_2d]
__global__ void Daxpy(stdmrg::complex<double> alpha,
                      const stdmrg::complex<double>* d_x,
                      stdmrg::complex<double>* d_y, int size) {
  unsigned int x = blockIdx.x * blockDim.x + threadIdx.x;
  unsigned int y = threadIdx.y + blockIdx.y * blockDim.y;
  int pitch = blockDim.x * gridDim.x;
  int idx = x + y * pitch;
  if ( idx < size ) {
    d_y[idx] = alpha * d_x[idx] + d_y[idx];
  }
}
\end{lstlisting}
The vector containing the final result of the addition overwrites the contents of the second vector operand $d_y$.

\begin{figure}[!ht]
\centering
  \begin{tikzpicture}[scale=0.75]
    \begin{axis}[
      height=8cm,
      width=10cm,
      xlabel=$size$,
      legend pos=north west,
      enlargelimits,
    ]
    \addplot[line width=2pt,dashed,black,mark=*] table[x index=0,y index=1,col sep=comma] {data/tikz/z_saxpy.txt};
    \addlegendentry{cpu time (ms)}
    \addplot[line width=2pt,red,mark=triangle] table[x index=0,y index=3,col sep=comma] {data/tikz/z_saxpy.txt};
    \addlegendentry{gpu\#0 time (ms)}
    \addplot[line width=2pt,blue,mark=*] table[x index=0,y index=5,col sep=comma] {data/tikz/z_saxpy.txt};
    \addlegendentry{gpu\#1 time (ms)}
    \end{axis}
  \end{tikzpicture}
  \begin{tikzpicture}[scale=0.75]
    \begin{axis}[
      height=8cm,
      width=10cm,
      xlabel=$size$,
      legend pos=south east,
      enlargelimits,
    ]
    \addplot[line width=2pt,dashed,black,mark=*] table[x index=0,y index=2,col sep=comma] {data/tikz/z_saxpy.txt};
    \addlegendentry{cpu GFlop/s}
    \addplot[line width=2pt,red,mark=triangle] table[x index=0,y index=4,col sep=comma] {data/tikz/z_saxpy.txt};
    \addlegendentry{gpu\#0 GFlop/s}
    \addplot[line width=2pt,blue,mark=*] table[x index=0,y index=6,col sep=comma] {data/tikz/z_saxpy.txt};
    \addlegendentry{gpu\#1 GFlop/s}
    \end{axis}
  \end{tikzpicture}

\caption{ZAXPY [left: time in ms, right: GFlops]}
\label{fig:tikz:z_saxpy}
\end{figure}
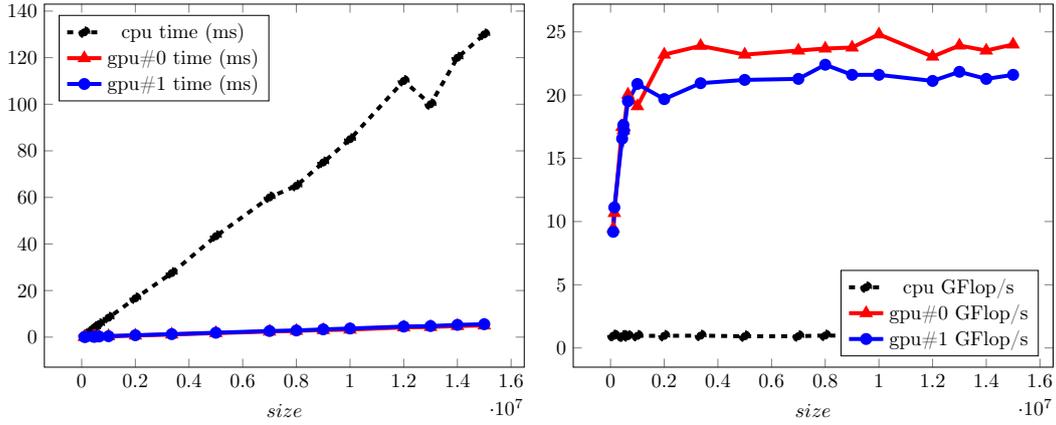

Table~\ref{tab:addition_vectors_tab} and Figure~\ref{fig:tikz:z_saxpy} report the complex number arithmetics with double precision execution time of our implementation for the \emph{Zaxpy} operation.
\begin{table}[!ht]  
\centering    
\renewcommand{\arraystretch}{1.3}
\renewcommand{\tabcolsep}{0.06cm}
\begin{tabular}{ccccccccc}
\hlinewd{1.0pt}
{\bf h} & {\bf cpu} & {\bf cpu} & {\bf gpu\#0} & {\bf gpu\#0} & {\bf gpu\#1} & {\bf gpu\#1} & {\bf ratio\#0} & {\bf ratio\#1} \\
{} & {\it time (ms)} & {\it Gflops} & {\it time (ms)} & {\it Gflops} & {\it time (ms)} & {\it Gflops} & {\it cpu/\#0} & {\it cpu/\#1} \\
\hlinewd{1.0pt}
100,000 & 0.83 & 0.97 & 0.08 & 9.46 & 0.09 & 9.18 & {\bf 9.77} & {\bf 9.49} \\
500,000 & 4.17 & 0.96 & 0.23 & 17.48 & 0.23 & 17.16 & {\bf 18.21} & {\bf 17.88} \\
1,000,000 & 8.33 & 0.96 & 0.42 & 19.12 & 0.38 & 20.88 & {\bf 19.92} & {\bf 21.75} \\
8,000,000 & 65.00 & 0.98 & 2.70 & 23.68 & 2.86 & 22.40 & {\bf 24.05} & {\bf 22.75} \\
10,000,000 & 85.00 & 0.94 & 3.23 & 24.80 & 3.70 & 21.60 & {\bf 26.35} & {\bf 22.95} \\
15,000,000 & 130.00 & 0.92 & 5.00 & 24.00 & 5.56 & 21.60 & {\bf 26.00} & {\bf 23.40} \\
\hlinewd{1.0pt}
\end{tabular}
\caption{Complex double precision Addition of vectors (ZAXPY)}
\label{tab:addition_vectors_tab}
\end{table}

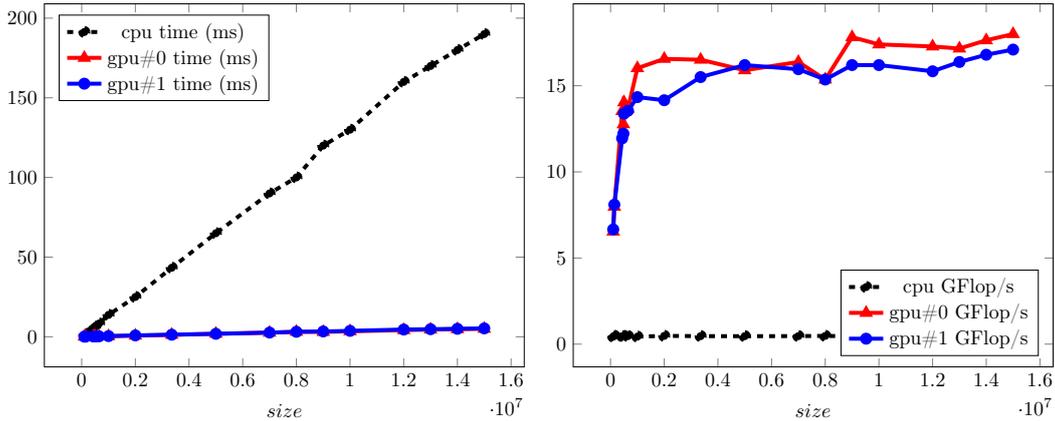
\begin{figure}[!ht]
\centering
  \begin{tikzpicture}[scale=0.75]
    \begin{axis}[
      height=8cm,
      width=10cm,
      xlabel=$size$,
      legend pos=north west,
      enlargelimits,
    ]
    \addplot[line width=2pt,dashed,black,mark=*] table[x index=0,y index=1,col sep=comma] {data/tikz/z_ewproduct.txt};
    \addlegendentry{cpu time (ms)}
    \addplot[line width=2pt,red,mark=triangle] table[x index=0,y index=3,col sep=comma] {data/tikz/z_ewproduct.txt};
    \addlegendentry{gpu\#0 time (ms)}
    \addplot[line width=2pt,blue,mark=*] table[x index=0,y index=5,col sep=comma] {data/tikz/z_ewproduct.txt};
    \addlegendentry{gpu\#1 time (ms)}
    \end{axis}
  \end{tikzpicture}
  \begin{tikzpicture}[scale=0.75]
    \begin{axis}[
      height=8cm,
      width=10cm,
      xlabel=$size$,
      legend pos=south east,
      enlargelimits,
    ]
    \addplot[line width=2pt,dashed,black,mark=*] table[x index=0,y index=2,col sep=comma] {data/tikz/z_ewproduct.txt};
    \addlegendentry{cpu GFlop/s}
    \addplot[line width=2pt,red,mark=triangle] table[x index=0,y index=4,col sep=comma] {data/tikz/z_ewproduct.txt};
    \addlegendentry{gpu\#0 GFlop/s}
    \addplot[line width=2pt,blue,mark=*] table[x index=0,y index=6,col sep=comma] {data/tikz/z_ewproduct.txt};
    \addlegendentry{gpu\#1 GFlop/s}
    \end{axis}
  \end{tikzpicture}

\caption{ZAXMY [left: time in ms, right: GFlops]}
\label{fig:tikz:z_ewproduct}
\end{figure}

The {\em element wise product} or {\em element by element product}, i.e., $y[i] = x[i] \times y[i]$, is intuitively parallel placing it as an excellent candidate for application on GPU. The product result of the elements of both vectors overwrites the corresponding elements of the second vector by a simple GPU kernel, presented in the following:
\begin{lstlisting}[language=cuda,label=lst:xmy_2d]
__global__ void ElementWiseProduct(
        stdmrg::complex<double> alpha,
        const stdmrg::complex<double>* d_x,
        stdmrg::complex<double>* d_y, int size) {
  unsigned int x = blockIdx.x * blockDim.x + threadIdx.x;
  unsigned int y = threadIdx.y + blockIdx.y * blockDim.y;
  int pitch = blockDim.x * gridDim.x;
  int idx = x + y * pitch;
  if ( idx < size ) {
    d_y[idx] = d_x[idx] * d_y[idx];
  }
}
\end{lstlisting}
Table~\ref{tab:ewp_vectors_tab} and Figure~\ref{fig:tikz:z_ewproduct} shows the double precision execution time of our implementation for the \emph{element wise product} operation.
\begin{table}[!ht]  
\centering    
\renewcommand{\arraystretch}{1.3}
\renewcommand{\tabcolsep}{0.10cm}
\caption{Complex double precision Element wise product (ZAXMY)}
\label{tab:ewp_vectors_tab}
\begin{tabular}{ccccccccc}
\hlinewd{1.0pt}
{\bf h} & {\bf cpu} & {\bf cpu} & {\bf gpu\#0} & {\bf gpu\#0} & {\bf gpu\#1} & {\bf gpu\#1} & {\bf ratio\#0} & {\bf ratio\#1} \\
{} & {\it time (ms)} & {\it Gflops} & {\it time (ms)} & {\it Gflops} & {\it time (ms)} & {\it Gflops} & {\it cpu/\#0} & {\it cpu/\#1} \\
\hlinewd{1.0pt}
100,000 & 1.37 & 0.44 & 0.09 & 6.53 & 0.09 & 6.66 & {\bf 14.92} & {\bf 15.21} \\
500,000 & 6.25 & 0.48 & 0.21 & 14.04 & 0.22 & 13.38 & {\bf 29.25} & {\bf 27.88} \\
1,000,000 & 13.75 & 0.44 & 0.37 & 16.02 & 0.42 & 14.34 & {\bf 36.71} & {\bf 32.86} \\
8,000,000 & 100.00 & 0.48 & 3.13 & 15.36 & 3.13 & 15.36 & {\bf 32.00} & {\bf 32.00} \\
10,000,000 & 130.00 & 0.46 & 3.45 & 17.40 & 3.70 & 16.20 & {\bf 37.70} & {\bf 35.10} \\
15,000,000 & 190.00 & 0.47 & 5.00 & 18.00 & 5.26 & 17.10 & {\bf 38.00} & {\bf 36.10} \\
\hlinewd{1.0pt}
\end{tabular}
\end{table}

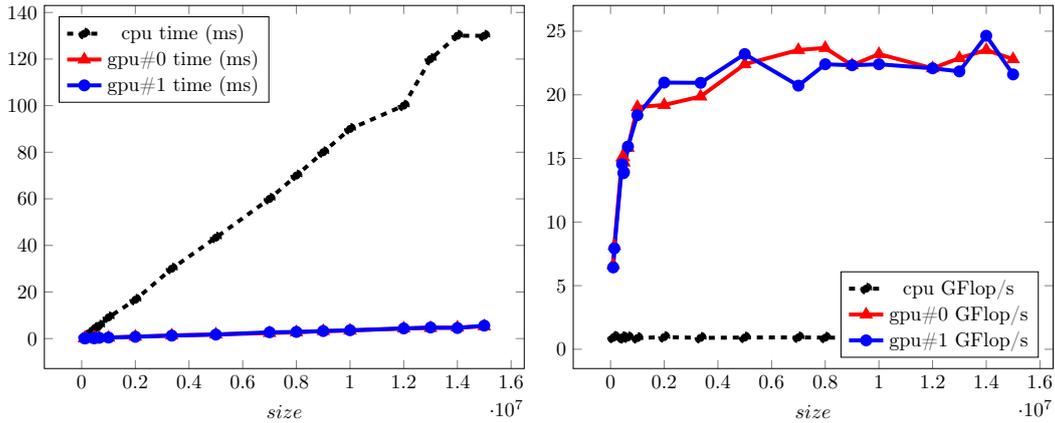
\begin{figure}[!ht]
\centering
  \begin{tikzpicture}[scale=0.75]
    \begin{axis}[
      height=8cm,
      width=10cm,
      xlabel=$size$,
      legend pos=north west,
      enlargelimits,
    ]
    \addplot[line width=2pt,dashed,black,mark=*] table[x index=0,y index=1,col sep=comma] {data/tikz/z_dot.txt};
    \addlegendentry{cpu time (ms)}
    \addplot[line width=2pt,red,mark=triangle] table[x index=0,y index=3,col sep=comma] {data/tikz/z_dot.txt};
    \addlegendentry{gpu\#0 time (ms)}
    \addplot[line width=2pt,blue,mark=*] table[x index=0,y index=5,col sep=comma] {data/tikz/z_dot.txt};
    \addlegendentry{gpu\#1 time (ms)}
    \end{axis}
  \end{tikzpicture}
  \begin{tikzpicture}[scale=0.75]
    \begin{axis}[
      height=8cm,
      width=10cm,
      xlabel=$size$,
      legend pos=south east,
      enlargelimits,
    ]
    \addplot[line width=2pt,dashed,black,mark=*] table[x index=0,y index=2,col sep=comma] {data/tikz/z_dot.txt};
    \addlegendentry{cpu GFlop/s}
    \addplot[line width=2pt,red,mark=triangle] table[x index=0,y index=4,col sep=comma] {data/tikz/z_dot.txt};
    \addlegendentry{gpu\#0 GFlop/s}
    \addplot[line width=2pt,blue,mark=*] table[x index=0,y index=6,col sep=comma] {data/tikz/z_dot.txt};
    \addlegendentry{gpu\#1 GFlop/s}
    \end{axis}
  \end{tikzpicture}

\caption{ZDOT [left: time in ms, right: GFlops]}
\label{fig:tikz:z_dot}
\end{figure}

One very expensive operation to do on CPUs for large vectors is the dot product. However, the most basic implementation of this operation, which uses a simple loop with simultaneous sums, is not very efficient on GPUs. This is why we split the dot product algorithm into two separate parts. The first one consists of the element by element parallel multiplication of the vectors, while the second one consists of the summation of all the results given by the first part. On a sequential processor, to implement the second part, a simple loop where a single variable is incremented is enough. Each element of the input data is handled by a thread, and the current sum (the partial result which is the sum of the $n^{th}$ first elements) is stored in the first thread of this block at the end of the reduction. To obtain the result of the dot product, the algorithm then returns the sum of all the partial sums of the different blocks.
In Table~\ref{tab:dot_product_tab} and Figure~\ref{fig:tikz:z_dot}, we compare the double precision execution time of our implementation for the \emph{dot product} on both CPU and GPU.
\begin{table}[!ht]  
\centering    
\renewcommand{\arraystretch}{1.3}
\renewcommand{\tabcolsep}{0.06cm}
\begin{tabular}{ccccccccc}
\hlinewd{1.0pt}
{\bf h} & {\bf cpu} & {\bf cpu} & {\bf gpu\#0} & {\bf gpu\#0} & {\bf gpu\#1} & {\bf gpu\#1} & {\bf ratio\#0} & {\bf ratio\#1} \\
{} & {\it time (ms)} & {\it Gflops} & {\it time (ms)} & {\it Gflops} & {\it time (ms)} & {\it Gflops} & {\it cpu/\#0} & {\it cpu/\#1} \\
\hlinewd{1.0pt}
100,000 & 0.88 & 0.91 & 0.12 & 6.63 & 0.12 & 6.43 & {\bf 7.27} & {\bf 7.05} \\
500,000 & 4.55 & 0.88 & 0.27 & 14.68 & 0.29 & 13.92 & {\bf 16.68} & {\bf 15.82} \\
1,000,000 & 9.09 & 0.88 & 0.42 & 19.04 & 0.43 & 18.40 & {\bf 21.64} & {\bf 20.91} \\
8,000,000 & 70.00 & 0.91 & 2.70 & 23.68 & 2.86 & 22.40 & {\bf 25.90} & {\bf 24.50} \\
10,000,000 & 90.00 & 0.89 & 3.45 & 23.20 & 3.57 & 22.40 & {\bf 26.10} & {\bf 25.20} \\
15,000,000 & 130.00 & 0.92 & 5.26 & 22.80 & 5.56 & 21.60 & {\bf 24.70} & {\bf 23.40} \\
\hlinewd{1.0pt}
\end{tabular}
\caption{Complex double precision Dot product (ZDOT)}
\label{tab:dot_product_tab}
\end{table}

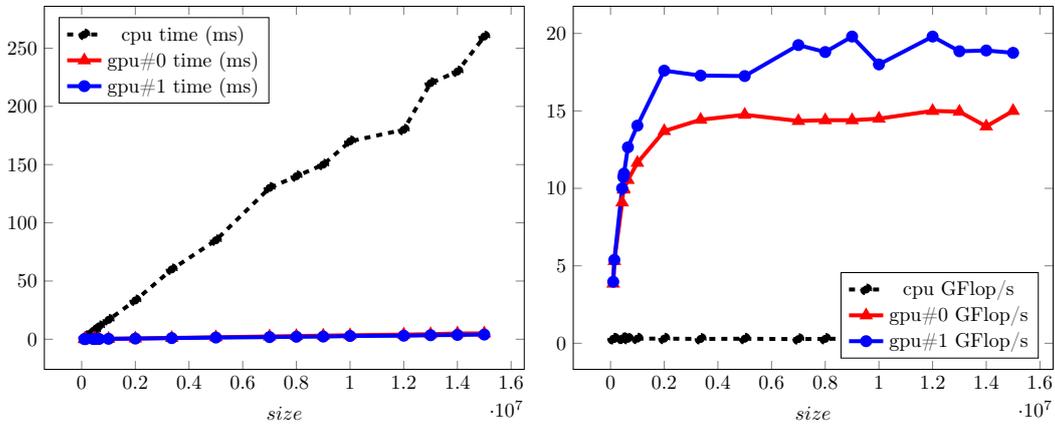
\begin{figure}[!ht]
\centering
  \begin{tikzpicture}[scale=0.75]
    \begin{axis}[
      height=8cm,
      width=10cm,
      xlabel=$size$,
      legend pos=north west,
      enlargelimits,
    ]
    \addplot[line width=2pt,dashed,black,mark=*] table[x index=0,y index=1,col sep=comma] {data/tikz/z_norm2.txt};
    \addlegendentry{cpu time (ms)}
    \addplot[line width=2pt,red,mark=triangle] table[x index=0,y index=3,col sep=comma] {data/tikz/z_norm2.txt};
    \addlegendentry{gpu\#0 time (ms)}
    \addplot[line width=2pt,blue,mark=*] table[x index=0,y index=5,col sep=comma] {data/tikz/z_norm2.txt};
    \addlegendentry{gpu\#1 time (ms)}
    \end{axis}
  \end{tikzpicture}
  \begin{tikzpicture}[scale=0.75]
    \begin{axis}[
      height=8cm,
      width=10cm,
      xlabel=$size$,
      legend pos=south east,
      enlargelimits,
    ]
    \addplot[line width=2pt,dashed,black,mark=*] table[x index=0,y index=2,col sep=comma] {data/tikz/z_norm2.txt};
    \addlegendentry{cpu GFlop/s}
    \addplot[line width=2pt,red,mark=triangle] table[x index=0,y index=4,col sep=comma] {data/tikz/z_norm2.txt};
    \addlegendentry{gpu\#0 GFlop/s}
    \addplot[line width=2pt,blue,mark=*] table[x index=0,y index=6,col sep=comma] {data/tikz/z_norm2.txt};
    \addlegendentry{gpu\#1 GFlop/s}
    \end{axis}
  \end{tikzpicture}

\caption{ZNORM [left: time in ms, right: GFlops]}
\label{fig:tikz:z_norm2}
\end{figure}

The results of the norm operation are given in Table~\ref{tab:normL2_tab} and Figure~\ref{fig:tikz:z_norm2}.
\begin{table}[!ht]  
\centering    
\renewcommand{\arraystretch}{1.3}
\renewcommand{\tabcolsep}{0.06cm}
\begin{tabular}{ccccccccc}
\hlinewd{1.0pt}
{\bf h} & {\bf cpu} & {\bf cpu} & {\bf gpu\#0} & {\bf gpu\#0} & {\bf gpu\#1} & {\bf gpu\#1} & {\bf ratio\#0} & {\bf ratio\#1} \\
{} & {\it time (ms)} & {\it Gflops} & {\it time (ms)} & {\it Gflops} & {\it time (ms)} & {\it Gflops} & {\it cpu/\#0} & {\it cpu/\#1} \\
\hlinewd{1.0pt}
100,000 & 1.72 & 0.29 & 0.13 & 3.85 & 0.13 & 3.97 & {\bf 13.28} & {\bf 13.69} \\
500,000 & 7.69 & 0.33 & 0.25 & 9.98 & 0.23 & 10.95 & {\bf 30.69} & {\bf 33.69} \\
1,000,000 & 16.67 & 0.30 & 0.43 & 11.65 & 0.36 & 14.05 & {\bf 38.83} & {\bf 46.83} \\
8,000,000 & 140.00 & 0.29 & 2.78 & 14.40 & 2.13 & 18.80 & {\bf 50.40} & {\bf 65.80} \\
10,000,000 & 170.00 & 0.29 & 3.45 & 14.50 & 2.78 & 18.00 & {\bf 49.30} & {\bf 61.20} \\
15,000,000 & 260.00 & 0.29 & 5.00 & 15.00 & 4.00 & 18.75 & {\bf 52.00} & {\bf 65.00} \\
\hlinewd{1.0pt}
\end{tabular}
\caption{Complex double precision NormL2 (ZNORM)}
\label{tab:normL2_tab}
\end{table}

As we can see in the presented results, GPU is clearly more effective than CPU with complex number arithmetics in double precision. In addition, for the operations considered on two different GPUs, the results are slightly similar, even if the first GPU is a little bit more efficient than the second. In the following, all experiments are performed on the second device, i.e., gpu\#0.

\section{Advanced Linear Algebra Operations}
\label{sec:advanced_linear_algebra_operations}

Many methods, such as finite element analysis, require handling widly large size sparse matrices, i.e., only few elements are non-zero values. To store efficiently these matrices on GPU memory, different structures exist, e.g., Compressed-sparse Row format. These structures are very important for Sparse matrix-vector product (SpMV) operation for instance, which is one of the most time consuming operation in sparse matrix computation. As explained before, the matrix-vector product proposed in this paper uses advanced auto-tuned techniques to organize threads on the CUDA grid. References~\cite{GPU:KK:2011,GPU:ZR:2012,GPU:ZR:2010,GPU:KGWHBA:2011,GPU:OSV:2010,GPU:DFG:2010} clearly demonstrated the efficiency of SpMV on GPU compared to CPU for real number arithmetics. Complex number arithmetics with double precision still remains a challenge.

We now reports in Table~\ref{tab:csr_spmv_tab} the SpMV execution time and the number of floating operations per second when using the CSR format within \emph{Alinea} for complex number arithmetics with double precision.
\begin{table}[!ht]  
\centering    
\renewcommand{\arraystretch}{1.3}
\renewcommand{\tabcolsep}{0.07cm}
\begin{tabular}{lcccccccc}
\hlinewd{1.0pt}
{\bf problem} & {\bf cpu} & {\bf cpu} & {\bf gpu\#0} & {\bf gpu\#0} & {\bf gpu\#1} & {\bf gpu\#1} & {\bf ratio\#0} & {\bf ratio\#1} \\
{} & {\it time (ms)} & {\it Gflops} & {\it time (ms)} & {\it Gflops} & {\it time (ms)} & {\it Gflops} & {\it cpu/\#0} & {\it cpu/\#1} \\
\hlinewd{1.0pt}
Audi3D-0 & 0.01 & 0.61 & 0.07 & 0.12 & 0.06 & 0.13 & {\bf 0.19} & {\bf 0.21} \\
Audi3D-1 & 0.20 & 0.67 & 0.11 & 1.23 & 0.12 & 1.07 & {\bf 1.84} & {\bf 1.60} \\
Audi3D-2 & 2.22 & 0.68 & 0.37 & 4.03 & 0.42 & 3.56 & {\bf 5.93} & {\bf 5.24} \\
Audi3D-3 & 20.00 & 0.71 & 2.22 & 6.41 & 3.03 & 4.70 & {\bf 9.00} & {\bf 6.60} \\
Audi3D-4 & 180.00 & 0.69 & 18.33 & 6.74 & 24.00 & 5.15 & {\bf 9.82} & {\bf 7.50} \\
Twingo3D-0 & 1.67 & 0.69 & 0.28 & 4.06 & 0.33 & 3.45 & {\bf 5.88} & {\bf 5.00} \\
Twingo3D-1 & 15.71 & 0.69 & 1.79 & 6.05 & 2.44 & 4.43 & {\bf 8.80} & {\bf 6.44} \\
Twingo3D-2 & 140.00 & 0.66 & 14.29 & 6.51 & 16.67 & 5.58 & {\bf 9.80} & {\bf 8.40} \\
Cylinder3D-0 & 0.37 & 0.67 & 0.13 & 1.90 & 0.15 & 1.67 & {\bf 2.84} & {\bf 2.51} \\
Cylinder3D-1 & 3.70 & 0.74 & 0.53 & 5.20 & 0.69 & 3.96 & {\bf 7.00} & {\bf 5.33} \\
Cylinder3D-2 & 36.67 & 0.69 & 4.00 & 6.30 & 5.26 & 4.79 & {\bf 9.17} & {\bf 6.97} \\
\hlinewd{1.0pt}
\end{tabular}
\caption{Double precision CSR matrix-vector multiplication}
\label{tab:csr_spmv_tab}
\end{table}
Numerical experiments clearly show that GPU operations are efficient than CPU operations.

\section{Iterative Krylov methods}
\label{sec:matrix_vector_product}

Knowing the efficiency of GPU architecture to perform linear algebra operations for complex number arithmetics with double precision, we now extend our analysis to Krylov iterative methods~\cite{GPU:LS:2010,GPU:KK:2011,GPU:BCK:2011,GPU:ZR:2012}. We have thus aimed to develop a preconditionned bi-conjugate gradient stabilized method (P-Bi-CGSTAB), a preconditionned P-BiCGSTAB parametered (l) and a preconditionned transpose-free quasi-minimal residual method (P-tfQMR)~\cite{saad_iterative_2003}, with optimized CUDA and dynamic auto-tuning on GPU.
Reference~\cite{cheikahamed:2012:inproceedings-2} exhibits the effectiveness of our template implementation for real number arithmetics compared to Cusp~\cite{GPU:CUSP:2010}, CUBLAS~\cite{GPU:CUBLAS}, CUSPARSE~\cite{GPU:CUSPARSE4.0:2011}.

In general, the transfer of data between the host (CPU) and the device (GPU), is one of the most time consuming aspects of the algorithm~\cite{GPU:CSVGM:2014}. In our Krylov algorithms implementation, all inputs are sent once from the CPU to the GPU before the iterative part starts. However, at each iteration, there has to be more than one dot product (or norm), which means that data must be copied from GPU to CPU.
The Krylov algorithm we have developed is launched on the host (CPU) but all the computing steps such as Zdot, Znorm, Zaxpy, or SpMV take place on the device (GPU). The implementations both on CPU and GPU are rigorously the same. All the results of the iterative Krylov methods presented in this paper are achieved with a residual tolerance threshold of $1\times10^{-9}$, an initial guess of zero and a maximum number of iterations equals to $1000$ for all of the presented methods. The aim of this experiments consists in comparing tha analysis both on CPU and GPU, with the same code, for complex number arithmetics with double precision.
Table~\ref{speed-up_Audi3D} shows the speed-up obtained for Audi3D.
\begin{table}[htbp]
\setlength\columnsep{0.1pt}
\centering
\begin{tabular}{lcccc}
\hlinewd{1.0pt}
problem & \#iter & CPU time (s) & GPU time (s) & speed-up \\
\hline
\multicolumn{3}{l}{\emph{P-BiCGSTAB}} & & \\
Audi3D-1 & 21 & 0.01 & 0.030 & \textbf{0.33} \\
Audi3D-2 & 53 &  0.24 & 0.106 & \textbf{2.26} \\
Audi3D-3 & 94 &  4.01 & 0.703 & \textbf{5.71} \\
Audi3D-4 & 183 & 85.70 & 9.209 & \textbf{9.31} \\
\hline
\multicolumn{3}{l}{\emph{P-BiCGSTAB(8)}} & & \\
Audi3D-1 & 6 & 0.03 & 0.110 & \textbf{0.27} \\
Audi3D-2 & 12 & 0.52 & 0.286 & \textbf{1.82} \\
Audi3D-3 & 31 & 12.47 & 2.162 & \textbf{5.77} \\
Audi3D-4 & 70 & 266.26 & 30.100 & \textbf{8.85} \\
\hline
\multicolumn{3}{l}{\emph{P-TFQMR}} & & \\
Audi3D-1 & 24 & 0.02 & 0.040 & \textbf{0.50} \\
Audi3D-2 & 52 & 0.27 & 0.113 & \textbf{2.40} \\
Audi3D-3 & 99 &  4.71 & 0.755 & \textbf{6.24} \\
Audi3D-4 & 214 &  102.17 & 10.786 & \textbf{9.47} \\
\hlinewd{1.0pt}
\end{tabular}
\caption{Speed-up of Audi3D}
\label{speed-up_Audi3D}
\end{table}
In Table~\ref{speed-up_Twingo3D} the speed-up for Twingo3D are collected. 
\begin{table}[htbp]
\setlength\columnsep{0.1pt}
\centering
\begin{tabular}{lcccc}
\hlinewd{1.0pt}
problem & \#iter & CPU time (s) & GPU time (s) & speed-up \\
\hline
\multicolumn{3}{l}{\emph{P-BiCGSTAB}} & & \\
Twingo3D-0 & 563 & 1.85 & 1.008 & \textbf{1.84} \\
Twingo3D-1 & 1000 & 29.45 & 5.730 & \textbf{5.14} \\
Twingo3D-2 & 1000 & 295.66 & 37.670 & \textbf{7.85} \\
\hline
\multicolumn{3}{l}{\emph{P-BiCGSTAB(8)}} & & \\
Twingo3D-0 & 1000 & 31.2 & 20.970 & \textbf{1.49} \\
Twingo3D-1 & 1000 & 273.81 & 54.630 & \textbf{5.01} \\
Twingo3D-2 & 1000 & 2559.67 & 324.500 & \textbf{7.89} \\
\hline
\multicolumn{3}{l}{\emph{P-TFQMR}} & & \\
Twingo3D-0 & 366 & 1.34 & 0.626 & \textbf{2.14} \\
Twingo3D-1 & 954 & 30.4 & 5.438 & \textbf{5.59} \\
Twingo3D-2 & 1000 & 318.93 & 38.090 & \textbf{8.37} \\
\hlinewd{1.0pt}
\end{tabular}
\caption{Speed-up of Twingo3D}
\label{speed-up_Twingo3D}
\end{table}
The CPU and GPU times in second, and the corresponding ratio are reported in Table~\ref{speed-up_Cylinder3D} 
\begin{table}[htbp]
\setlength\columnsep{0.1pt}
\centering
\begin{tabular}{lcccc}
\hlinewd{1.0pt}
problem & \#iter & CPU time (s) & GPU time (s) & speed-up \\
\hline
\multicolumn{3}{l}{\emph{P-BiCGSTAB}} & & \\
Cylinder3D-0 & 49 & 0.040 & 0.060 & \textbf{0.67} \\
Cylinder3D-1 & 86 & 0.680 & 0.216 & \textbf{3.15} \\
Cylinder3D-2 & 162 & 12.290 & 1.886 & \textbf{6.52} \\
\hline
\multicolumn{3}{l}{\emph{P-BiCGSTAB(8)}} & & \\
Cylinder3D-0 & 10 & 0.090 & 0.182 & \textbf{0.50} \\
Cylinder3D-1 & 19 & 1.450 & 0.521 & \textbf{2.79} \\
Cylinder3D-2 & 71 & 47.670 & 7.460 & \textbf{6.39} \\
\hline
\multicolumn{3}{l}{\emph{P-QMR}} & & \\
Cylinder3D-0 & 85 & 0.090 & 0.106 & \textbf{0.85} \\
Cylinder3D-1 & 196 & 2.050 & 0.482 & \textbf{4.25} \\
Cylinder3D-2 & 282 & 23.870 & 3.302 & \textbf{7.23} \\
\hlinewd{1.0pt}
\end{tabular}
\caption{Speed-up of Cylinder3D}
\label{speed-up_Cylinder3D}
\end{table}
The numerical results of the considered iterative Krylov methods confirm the efficiency of GPU computations compared to CPU for solving sparse linear systems. In addition, the ratio increases when the size of the problem increases for all cases, i.e., when the mesh size is more fine, GPU is more efficient compared to CPU. Unfortunaltely, when the mesh is too fine, the corresponding assembled matrix becomes too large for GPU memory\ldots In this case, domain decomposition method~\cite{SBG1996},~\cite{quarteroni_domain_1999},~\cite{toselli_domain_2004},~\cite{NME:NME1620320604},~\cite{magoules:journal-auth:21},~\cite{magoules:journal-auth:16} based on iterative methods is an issue.
The Schwarz method~\cite{magoulesf_contrib_3:Lions:1988:SAM},~\cite{magoulesf_contrib_3:Lions:1989:SAM},~\cite{magoulesf_contrib_3:Lions:1990:SAM},~\cite{cai_overlapping_1998} have encountered a success for solving large size problem.
To faster the convergence, many references~\cite{magoulesf_contrib_3:Chevalier:1998:SMO,magoulesf_contrib_3:Gander:2000:OSM,magoules:journal-auth:23,magoules:journal-auth:18,magoules:journal-auth:14} exhibit the interest of optimizing the interface conditions between the subdomains.
In order to use this approach for acoustic problems modeled by the Helmholtz equation, continuous optimized interface conditions between the subdomains must be developped as in~\cite{magoules:journal-auth:9},~\cite{magoules:journal-auth:10},~\cite{magoules:journal-auth:13},~\cite{magoules:journal-auth:28}. 
Alternative discrete optimization techniques as introduced in~\cite{magoules:proceedings-auth:6}, \cite{magoules:journal-auth:29},~\cite{magoules:journal-auth:12},~\cite{magoules:journal-auth:30},~\cite{magoules:journal-auth:20},~\cite{magoules:journal-auth:17}.

Authors present how domain decomposition method is designed efficiently on GPU in~\cite{cheikahamed:2013:inproceedings-4,ahamed2013stochastic} and proved in~\cite{ahamed2013stochastic} the interest of Schwarz methods on a cluster of GPUs.

\section{Concluding remarks}
\label{sec:concluding_remarks}

This paper presented the performance evaluation and analysis of linear
algebra operations together with their uses within Krylov methods for solving acoustic problem on Graphics Processing Unit (GPU) for complex number arithmetics with double precision. The numerical experiments have been performed in two different system accelerated generations of nVidia graphics card: GTX570 and Tesla K20c.
We have used matrices arising from the finite element modeling of the acoustic within a cylinder and a car compartment. The presented results clearly demonstrate the interest of the use of GPU device to compute linear algebra operations, and outline the robustness, performance
and efficiency of solve the Helmholtz equations for acoustic problems.

\bibliography{bib/hpcc2014_paper2_ac_fm,bib/MAGOULES-JOURNAL1,bib/MAGOULES-PROCEEDINGS1}
\bibliographystyle{abbrv}

\end{document}